\documentclass[12pt]{article} \textwidth=6in \oddsidemargin=0in
\usepackage{amssymb} 
\begin{document}\def\ov{\over} \def\x{\xi} \def\t{\tau}
\def\inv{^{-1}} \def\cd{\cdots} \def\ld{\ldots} \def\r{\rho}
\def\({\Bigg(} \def\){\Bigg)}\def\lb{\Big[} \def\rb{\Big]} 
\def\iy{\infty} \def\be{\begin{equation}}\def\ee{\end{equation}} 
\def\s{\sigma} \def\bs{\backslash} \def\eq{\equiv}
\newcommand{\twotwo}[4]{\left(\begin{array}{cc}#1&#2\\&\\#3&#4\end{array}\right)} \newcommand{\onetwo}[2]{\left(\begin{array}{c}#1\\\\#2\end{array}\right)} \def\nb{\bar n} \def\Z{\mathbb Z}
\def\l{\ell}\def\ph{\varphi} \def\m{\mu} \def\n{\nu} \def\a{\alpha}
\def\mone{\left(\begin{array}{c}1\\1\\\vdots\end{array}\right)} \def\noi{\noindent} \def\P{\mathbb P} \def\t{\tau} \newcommand{\br}[2]{\left[{#1\atop #2}\right]_\t} \def\ep{\varepsilon} \newcommand{\C}[1]{{\cal C}_{#1}} \def\onem{(1\ 0\ 0\ \cd)} \def\A{\mathcal{A}}
\def\B{\mathcal{B}} \def\F{\mathcal{F}} 
\hfill September 21, 2010

\begin{center}{\large \bf On ASEP with Periodic Step Bernoulli Initial Condition}\end{center}

\begin{center}{\large\bf Craig A.~Tracy}\\
{\it Department of Mathematics \\
University of California\\
Davis, CA 95616, USA\\
email: tracy@math.ucdavis.edu}\end{center}

\begin{center}{\large \bf Harold Widom}\\
{\it Department of Mathematics\\
University of California\\
Santa Cruz, CA 95064, USA\\
email: widom@ucsc.edu}\end{center}

\begin{center}{\bf I. Introduction}\end{center}

In the asymmetric simple exclusion process (ASEP) on the integers $\Z$ a particle waits exponential time, then moves to the right with probability $p$ if that site is unoccupied (or else stays put) or to the left with probability $q=1-p$ if that site is unoccupied (or else stays put). A formula was derived in \cite{TW1} for $\P_Y(x_\l(t)\le x)$, the probability distribution function for $x_\l(t)$, the position of the $\l$th particle from the left at time~$t$, given the initial configuration $Y=\{y_1<y_2<\cd\}$. The set $Y$ may be finite or infinite on the right. The formula will be stated below. It is a sum over all (finite) subsets $S$ of $Y$ of integrals of order $|S|$. 

For step initial condition, when $Y=\Z^+$ (the positive integers), three pleasant things happen:

(i) For each $k$ one can sum over all $S$ with $|S|=k$, thus replacing a sum over all finite $S\subset\Z^+$ by a sum over $k$. 

(ii) There is a combinatorial identity that replaces the integrands by simpler ones which are symmetric in the integration variables.

(iii) The resulting $k$-dimensional integrals are coefficients in the expansion of a Fredholm determinant.

This led to a representation of the probability distribution in terms of a Fredholm determinant \cite{TW2}, which was amenable to asymptotic analysis \cite{TW3,SS1,SS2,ACQ}.  

With {\it step Bernoulli initial condition} $Y$ is not deterministic, but each site in $\Z^+$ is, independently of the others, initially occupied with probability $\r$. In this case also we have (i)--(iii), and this made asymptotic analysis possible \cite{TW4,CQ}.

In another direction, Lee \cite{L} considered the cases when $Y=m\Z^+$, and found that (i) and (ii) hold when $m=2$, but only (i) holds when $m>2$. Even for $m=2$, (iii) does not seem to hold since the $k$-dimensional integrals are not, or not obviously, related to a determinant or Pfaffian.\footnote{Earlier work on TASEP with initial condition $2\Z^+$ or $2\Z$ is in \cite{S,BFPS,BFS}.} 

Nevertheless it seems worthwhile to consider a generalization, which we call\linebreak
{\it periodic step Bernoulli initial condition}. Here we have an $m$-periodic function $n\to\r_n$ from $\Z$ to $[0,\,1]$ and assume that initially a site $n\in\Z^+$ is occupied with probability $\r_n$. The most that can be expected in this generality is that (i) holds, and we shall show that it does.

Here are the formulas alluded to and the result of the present paper. We assume $q\ne0$, set $\t=p/q$, and recall that the $\t$-binomial coefficient $\br{N}{\l}$ is defined as
\[{(1-\t^N)\,(1-\t^{N-1})\cdots (1-\t^{N-\l+1})\ov (1-\t)\,(1-\t^2)\cdots (1-\t^\l)}.\]

We define
\[\ep(\x)=p\,\x\inv+q\,\x-1,\ \ \ \ f(\x_i,\,\x_j)={\x_j-\x_i\ov p+q\x_i\x_j-\x_i},\] 
and then
\[I(x,k,\x)=I(x,k,\x_1,\ld,\x_k)=\prod_{i<j}f(\x_i,\,\x_j)\;\prod_i{\x_i^{x}\,e^{\ep(\x_i)t}\ov 1-\x_i}.\]
All indices in the products run from 1 to $k$. Notice that $I(x,k,\x)$ depends on $t$, although it is not displayed in the notations.

Finally, given two sets of integers $U$ and $V$ we define
\[\s(U,\,V)=\#\{(u,\,v): u\in U,\ v\in V,\ {\rm and}\ u\ge v.\}.\]

Theorem 5.2 of \cite{TW1}, slightly modified, is the formula
\be\P_Y(x_\l(t)\le x)=\sum_{{S\subset Y\atop |S|\ge \l}}\,c_{\l,k}\,\t^{\s(S,\,Y)}\,\int_{\C{R}}\cd\int_{\C{R}} I(x,k,\x)\,\prod_{i=1}^k\x_i^{-s_i}\,d\x_1\cd d\x_k,\label{PY}\ee 
where $k=|S|$ and
\[S=\{s_1,\ld,s_k\},\ \ \ \ \ (s_1<s_2<\cd<s_k)\]
\[ c_{\l,k}=(-1)^{\l}\,q^{k(k-1)/2}\,\t^{\l(\l-1)/2-k\l}\,
\br{k-1}{\l-1}.\]
The sum is taken over all (finite) subsets $S$ of $Y$ with $|S|\ge \l$. 
The contour $\C{R}$ is the circle with center 0 and radius $R$, which is assumed so large that the denominators $p+q\x_i\x_j-\x_i$  are nonzero on and outside the contours.  (All integrals are given the factor $1/2\pi i$.)
  
When $Y$ is not deterministic the probability $\P(x_\l(t)\le x)$ is obtained by averaging the right side of (\ref{PY}) over all $Y$. The focus of the present note is (i), the computation of the sum over all $S\subset Y$ with $|S|=k$, and then the average over all initial configurations $Y\subset\Z^+$. What will be important here are only those ingredients of (\ref{PY}) that depend on $Y$ and  $S$. These combine as
\be\t^{\s(S,\,Y)}\ \prod_{i=1}^k\x_i^{-s_i}.\label{tauxi}\ee

For step initial condition the sum of (\ref{tauxi}) over all $S$ equals \cite[p.~838]{TW1}
\[\t^{k(k+1)/2}\,\prod_{i=1}^k {1\ov \x_i\cd\x_k-\t^{k-i+1}}.\]
For step Bernoulli initial condition the sum over $S\subset Y$ followed by the average over initial configurations $Y$ is \cite[p.~830]{TW4}
\be\t^{k(k+1)/2}\,\prod_{i=1}^k {\r\ov \x_i\cd\x_k-1+\r-\t^{k-i+1} \r}.\label{step}\ee
For $Y=m\Z^+$ the sum of (\ref{tauxi}) over $S$ equals \cite{L}
\be\t^{k(k+1)/2}\prod_{i=1}^k{1\ov (\x_i\cd\x_k)^m-\t^{k-i+1}}.\label{alternate}\ee

To state the formula we derive here for periodic step Bernoulli initial condition we define 
\[\ph(i,n)=1-\r_n+\r_n\,\t^{k-i+1}\]
and $m\times m$ matrices $\A_i$ with entries 
\be A_{i,\,\m,\,\n}={1\ov(\x_i\cd\x_k)^{m}-\prod\limits_{n=1}^m\ph(i,\,n)}
\times\left\{\begin{array}{ll}(\x_i\cd\x_k)^{m}\,\x_i^{-\n}\,\r_{\n}\,
\prod\limits_{n=\m+1}^{\n-1}\ph(i,\,n)&{\rm if}\ \m<\n,\\&\\ \x_i^{-\n}\,\r_{\n}\,
\prod\limits_{n=\m+1}^{\n+m-1}\ph(i,\,n)&{\rm if}\ \m\ge\n.
\end{array}\right.\label{A}\ee
The indices $\m$ and $\n$ run from 0 to $m-1$. 

The result is that the average over $Y$ of the sum of (\ref{tauxi}) over $S\subset Y$ with $|S|=k$ is equal to $\t^{k(k+1)/2}$ times
\be(1\ 0\ \cd\ 0)\ \prod_{i=1}^k\,\A_i\ \left(\begin{array}{c}1\\ 
\vdots\\1\end{array}\right).\label{result}\ee
This is, in words, the sum of the entries in the top row of the matrix product. In the product, matrices with lower index $i$ are on the left.
To recapitulate:

\noi{\bf Theorem}. For periodic step Bernoulli initial condition the probability $\P(x_\l(t)\le x)$ is given by the right side of (\ref{PY}) with $\sum_{S\subset Y}$ replaced by $\sum_{k\ge \l}$ and the product (\ref{tauxi}) replaced by $\t^{k(k+1)/2}$ times (\ref{result}).

In the simplest special cases, $\r_n=\r$ when $n\equiv\n$ (mod $m$) and $\r_n=0$ otherwise. If we take $0\le\n<m$, then (\ref{result}) is equal to
\[A_{1,0,\n}\,\prod_{i=2}^k A_{i,\n,\n},\]
which we find is equal to
\[\prod_{i=1}^k{\r\ov(\x_i\cd\x_k)^m-(1-\r+\r\,\t^{k-i+1})}\]
when $\n=0$ and
\[(\x_1\cd\x_k)^{m-\n}\,\prod_{i=1}^k{\r\ov(\x_i\cd\x_k)^m-(1-\r+\r\,\t^{k-i+1})}\]
otherwise. Alternatively, it is given by just the last formula if instead we take\linebreak $0<\n\le m$.

In the course of the proof we obtain an analogous statement for what we may call {\it general step Bernoulli initial condition} in which the function $\r:\Z^+\to[0,\,1]$ is arbitrary. We define $\B_i$ to be the semi-infinite matrix with entries
\be B_{i,\,\m,\,\n}=\left\{\begin{array}{ll}\x_i^{-\n}\,\r_{\n}\,
\prod\limits_{n=\m+1}^{\n-1}\ph(i,\,n)&{\rm if}\ \m<\n,\\&\\0&{\rm if}\ \m\ge\n,
\end{array}\right.\label{B}\ee 
the indices staisfying $0\le\m,\,\n<\iy$.
We shall see that the average over $Y$ of the sum of (\ref{tauxi}) over $S\subset Y$ with $|S|=k$ is then equal to $\t^{k(k+1)/2}$ times
\be\onem\ \prod_{i=1}^k\,\B_i\ \mone.\label{general}\ee

Observe that this is what is obtained by formally letting $m\to\iy$ in the periodic case. Observe also that if $\r$ takes only the values 0 or 1 we are back in the case of deterministic initial configuration, and what we obtain can be seen to be just a restatement of (\ref{PY}).
 
\begin{center}{\bf II. Proof of the theorem}\end{center}

We assume first that $\r:\Z^+\to[0,\,1]$ is arbitrary and want to compute the average over $Y$ of the sum of (\ref{tauxi}) over $S\subset Y$ with $|S|=k$. As in \cite{TW4} we do this in the opposite order. We first take a fixed $S$ and average over all $Y\supset S$. (Later we will take the sum over $S$.) Also, as in \cite{TW4}, we may assume that $Y\subset[1,\,N]$, and at the end  we let $N\to\iy$. (In fact all we shall use is that $N\ge s_k$.)

\begin{center}{\it 1. The average over Y}\end{center}

The part of (\ref{tauxi}) that depends on $Y$ is $\t^{\s(S,\,Y)}$, while the probability of an initial configuration $Y$ is
\[\prod_{n\in Y}\r_n\;\prod_{n\in [1,\,N]\bs Y}(1-\r_n),\]
so we have to compute
\[\sum_{S\subset Y\subset [1,\,N]}\,\t^{\s(S,\,Y)}\,
\prod_{n\in Y}\r_n\;\prod_{n\in [1,\,N]\bs Y}(1-\r_n).\]
Any element of $Y$ which is larger than $s_k$ does not affect $\s(S,\,Y)$, so this sum may be written as a product of sums,
\be\(\sum_{S\subset Y\subset [1,\,s_k]}\t^{\s(S,\,Y)}\,
\prod_{n\in Y}\r_n\;\prod_{n\in [1,\,s_k]\bs Y}(1-\r_n)\)\cdot
\(\sum_{Y\subset(s_k,\,N]}\prod_{n\in Y}\r_n\;\prod_{n\in (s_k,\,N]\bs Y}(1-\r_n)\).\label{sumprod}\ee

The second factor equals 1, since it is  
\[\prod_{n\in(s_k,\,N]}(\r_n+(1-\r_n)),\] 
when written as a sum of products.

It remains to consider the first factor in (\ref{sumprod}).  
We use the notation
\[Y_i=(s_{i-1},\,s_i)\cap Y,\ \ \ (s_0=0),\]
so $Y$ is the disjoint union $S\cup Y_1\cup\cd\cup Y_k$. Since the number of elements of $S$ greater than or equal to $s_i$ is $k-i+1$, we have
\[\s(S,\,Y)=\s(S,\,S)+\sum_{i=1}^k\s(S,\,Y_i)=k(k+1)/2+\sum_{i=1}^k (k-i+1)\,|Y_i|.\]
Therefore the first factor in (\ref{sumprod}) equals 
\[\t^{k(k+1)/2}\;\prod_{n\in S}\r_n\  \prod_{i=1}^k\;\(\sum_{Y_i\subset (s_{i-1},\,s_i)}\t^{(k-i+1)\,|Y_i|}\,
\prod_{n\in Y_i}\r_n\;\prod_{n\in (s_{i-1},\,s_i)\bs Y_i}(1-\r_n)\)\]
\[=\t^{k(k+1)/2}\;\prod_{n\in S}\r_n\ \prod_{i=1}^k\;\(\sum_{Y_i\subset (s_{i-1},\,s_i)}\,
\prod_{n\in Y_i}(\t^{(k-i+1)}\r_n)\;\prod_{n\in (s_{i-1},\,s_i)\bs Y_i}(1-\r_n)\).\]

The sum over all $Y_i\subset (s_{i-1},\,s_i)$ is 
\[\prod_{n\in(s_{i-1},\,s_i)}\ph(i,n)=\prod_{n\in(s_{i-1},\,s_i)}(\t^{k-i+1}\,\r_n+(1-\r_n)),\]
as is seen by writing this as a sum of products.
Hence (\ref{sumprod}) equals $\t^{k(k+1)/2}$ times
\be\prod_{n\in S}\r_n\;\cdot\;\prod_{i=1}^k\;\prod_{n\in(s_{i-1},\,s_i)}\ph(i,n).\label{lucky}\ee

\begin{center}{\it 2. The sum over S: general $\r$}\end{center}

This gives the average of $\t^{\s(S,\,Y)}$ over all $Y\supset S$. Now we are to multiply this by $\prod_i\,\x_i^{-s_i}$, the second factor in (\ref{tauxi}), and sum over all $S$ with $|S|=k$. Since $S=\{s_1,\ld,s_k\}$, the sum over all $S$ equals
\be\sum_{0<s_1<\cd<s_k}\ \prod_{i=1}^k\;\(\r_{s_i}\;\x_i^{-s_i}\prod_{n\in(s_{i-1},\,s_i)}\ph(i,n)\).\label{Ssum}\ee
In terms of the entries of the matrices $\B_i$ given by (\ref{B}) this equals
\[\sum_{s_1,\ld, s_k=0}^\iy B_{1,\,0,\,s_1}\,B_{2,\,s_1,\,s_2}\cd B_{k,\,s_{k-1},\,s_k},\]
which equals (\ref{general}), the sum of the entries of the top row of the matrix product. (We used here that $B_{i,\,\m,\,\n}=0$ unless $\m<\n$.)

\begin{center}{\it 3. The sum over S: periodic $\r$}\end{center}

We shall rewrite the sum (\ref{Ssum}) using the periodicity of $\r$. There are unique representations
\[s_i=u_i\,m+v_i,\ \ \ {\rm where}\ \ u_i\ge0\ \ {\rm and}\ \ 0\le v_i<m.\]
We have $s_{i-1}<s_i$ if and only if the following hold:
\[{\rm if}\ \ v_{i-1}<v_i\ \ {\rm then}\ \ u_{i-1}\le u_i,\]
\[{\rm if}\ \ v_{i-1}\ge v_i\ \ {\rm then}\ \ u_{i-1}<u_i.\]
If we set $u_i=t_1+\cd+t_i$, so that
$s_i=(t_1+\cd+t_i)\,m+v_i$,
then the above becomes
\be\begin{array}{ll}t_i\ge0&{\rm if}\ v_{i-1}<v_i\\&\\ t_i>0&{\rm if} \ v_{i-1}\ge v_i.\end{array}\label{tv}\ee

With this notation the factor $\x_i^{-s_i}$ in (\ref{Ssum}) becomes
\[\x_i^{-(t_1+\cd+t_i)\,m-v_i},\]
the factor $\r_{s_i}$ becomes $\r_{v_i}$ by the periodicity of $\r$, and by the periodicity of $\ph(i,\,n)$ in~$n$ the last factor becomes
\[\prod_{n=v_{i-1}+1}^{t_i\,m+v_i-1}\,\ph(i,n)=\left\{\begin{array}{ll}\prod\limits_{n=v_{i-1}+1}^{v_i-1}\,\ph(i,n)\cdot
\Big(\prod\limits_{n=1}^{m}\,\ph(i,n)\Big)^{t_i}&{\rm if}\ v_{i-1}<v_i,\\&\\
\prod\limits_{n=v_{i-1}+1}^{v_i+m-1}\,\ph(i,n)\cdot\Big(\prod\limits_{n=1}^{m}\,\ph(i,n)\Big)^{t_i-1}&{\rm if}\ v_{i-1}\ge v_i.\end{array}\right.\]

If we fix the $v_i$ and first take the sum in (\ref{Ssum}) over the $t_i$ according to (\ref{tv}) we obtain
\[{1\ov (\x_i\cd\x_k)^{m}-\prod\limits_{n=1}^{m}\,\ph(i,n)}\times
\left\{\begin{array}{ll}(\x_i\cd\x_k)^{m}\,\xi_i^{-v_i}\,\r_{v_i}\,\prod\limits_{n=v_{i-1}+1}^{v_i-1}\,\ph(i,n)&{\rm if}\ v_{i-1}<v_i,\\
&\\ \xi_i^{-v_i}\,\r_{v_i}\,\prod\limits_{n=v_{i-1}+1}^{v_i+m-1}\,\ph(i,n)&{\rm if}\ v_{i-1}\ge v_i.\end{array}\right.\]

Then the product of this over $i$ is to be summed over all $v_i$ satisfying $0\le v_i<m$. In terms of the entries of the matrices $\A_i$ given by (\ref{A}) the sum equals
\[\sum_{v_1,\ld, v_k=0}^{m-1} A_{1,\,0,\,v_1}\,A_{2,\,v_1,\,v_2}\cd A_{k,\,v_{k-1},\,v_k},\]
which equals (\ref{result}). This completes the proof of the theorem.

\begin{center}{\bf Acknowledgment}\end{center}

This work was supported by the National Science Foundation through grants DMS-0906387 (first author) and DMS-0854934 (second author).

\end{document}